\documentclass[11pt]{article}
\usepackage{enumerate}
\usepackage{amsmath,epsfig,amssymb,amsbsy,verbatim,array,color,graphics,graphicx}
\usepackage{amssymb,amsmath,amsthm,amsfonts}
\usepackage{graphicx}
\usepackage{dsfont}

\newcount\mwhr\newcount\mwmin
\mwhr=\number\time\divide
\mwhr by60 \mwmin=\mwhr\multiply\mwmin by-60
\advance\mwmin by\number\time{}
\newcommand\R{\mathbb{R}}
\newcommand\Z{\mathbb{Z}}

\begin{document}

\title{\bf Sums of Two Squares Visualized}
\author{Ishan Banerjee\footnote{Chennai Mathematical Institute, Chennai, India. Email: \texttt{ishanbanerjee.314@gmail.com}.}
\and
Amites Sarkar\footnote{Western Washington University, Bellingham, USA. Email: \texttt{amites.sarkar@wwu.edu}.}}

\maketitle

\begin{abstract}
We provide a geometric interpretation of Brillhart's celebrated algorithm for expressing a prime $p\equiv 1\pmod 4$
as the sum of two squares.
\end{abstract}

\section{Introduction}

On Christmas Day 1640, Pierre de Fermat stated his famous ``two-square" theorem in a letter to Marin Mersenne.
This theorem, which was first stated by Girard in 1625, and first proved by Euler in 1749, states that an odd
prime $p$ is the sum of two squares if (and only if) $p\equiv 1\pmod 4$. Among the many admirers of the theorem
was G.H. Hardy, who wrote in {\it A Mathematician's Apology} that the theorem is one of the finest of arithmetic,
and also that ``there is no proof within the comprehension of anybody but a fairly expert mathematician."

\bigskip

\noindent Perhaps Hardy meant to encourage his younger readers, and flatter
his older ones! In any case, his famous textbook with E. M. Wright {\it An Introduction to the Theory of Numbers}
contains a neat geometrical proof of the two-square theorem, due to J.H. Grace~\cite{Grace}. Grace's proof can be
summarized as follows. First, we solve the congruence $u^2+1\equiv 0\pmod p$. (No major expertise required here:
writing $p=4k+1$, a short argument using Wilson's theorem shows that $u=(2k)!$ is a solution. The proof is due
to Lagrange.) Next consider the integer lattice $L$ given by
\begin{equation}\label{eq:Grace}
L=\{(x,y)\in\Z^2: x\equiv uy\pmod p\}.
\end{equation}
$L$ is a lattice containing a proportion $1/p$ of integer points; more precisely, the area of any fundamental
parallelogram of the lattice is $p$. But $L$ is also a {\it square lattice}, since if $P=(g,h)$ is a point
of $L$ closest to the origin $O$, then $g\equiv uh\pmod p$, so that $ug\equiv -h\pmod p$, and so $(-h,g)\in L$ too.
These two facts together imply that the area $g^2+h^2$ of the square on side $OP$ is $p$, proving the theorem.

\bigskip

\noindent All well and good. But is there a fast algorithm for finding $g$ and $h$? Indeed there is. In 1972,
John Brillhart~\cite{Brill}, building on earlier work of Legendre, Serret and Hermite, proposed the following method,
starting with $p$ and the number $u$ defined above. If $p=u^2+1$, stop. Otherwise, simply apply the Euclidean algorithm
to $p$ and $u$, and stop when the first remainder $r$ satisfies $r<\sqrt{p}$. If the {\it next} remainder is $r'$, then
$(g,h)=(r,r')$. (The fly in the ointment is that we first need to find $u$. Lagrange's formula $u=(2k)!$ is computationally
inefficient, but Shiu~\cite{Shiu} provides a non-deterministic algorithm with a short expected running time; the key is to
first find a quadratic non-residue modulo $p$, and then repeatedly square it.)

\bigskip

\noindent An example will serve to illustrate the method. Let $p=73$, and solve $u^2+1\equiv 0\pmod p$ to get $u=\pm27$.
Taking $u=27$, the Euclidean algorithm now yields:
\begin{align*}
73&=2\cdot 27+19\\
27&=1\cdot 19+{\bf 8}\\
19&=2\cdot 8+{\bf 3}
\end{align*}
and, lo and behold, $8^2+3^2=73$.

\bigskip

\noindent Magic? Maybe. Certainly, Brillhart's proof of his algorithm's correctness is somewhat mysterious,
even after simplification by Wagon~\cite{Wagon}. Shiu~\cite{Shiu} proves it using properties of continuants.
But a closer look at Grace's proof reveals that Brillhart's algorithm has a secret geometric twin: {\it Lagrange's
algorithm}.

\section{Lattice reduction}

\noindent Two integer vectors in $\R^2$, such as ${\bf b_1}=(3,1)^{\sf T}$ and ${\bf b_2}=(5,3)^{\sf T}$, generate a lattice $L$
in the plane, consisting of all vectors of the form $m{\bf b_1}+n{\bf b_2}$, where $m$ and $n$ are integers. Noting that
$||{\bf b_1}||=\sqrt{10}$ and $||{\bf b_2}||=\sqrt{34}$, we might wonder if there is a better way of generating the same
lattice. What do we mean by better? One answer is that we should try to minimize $||{\bf b_1}||$ and $||{\bf b_2}||$.
Equivalently, we need to find the two points of $L$ closest to the origin $O$.
It turns out that there is a simple way to do this, which was discovered independently by Lagrange (first) and Gauss (second).
We will call it {\it Lagrange's algorithm}. It is similar to the Gram-Schmidt process, but more complicated, since we are only
allowed to adjust vectors by {\it integer} multiples of other vectors, because we need to stay inside the lattice $L$. Indeed,
given any two linearly independent vectors in $\R^2$, the Gram-Schmidt process terminates after one step, while Lagrange's
algorithm will terminate, but, depending on the initial two vectors, could take arbitrarily many steps to do so.

\bigskip

\begin{figure}[htp!]
\centering
\includegraphics[width=0.6\linewidth]{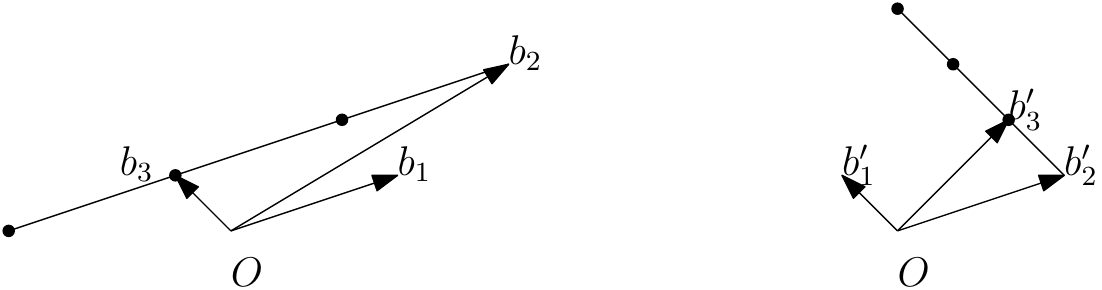}
\caption{\small Two steps of Lagrange's algorithm}
\label{2step}
\end{figure}

\noindent Figure~\ref{2step} illustrates how Lagrange's algorithm works. We start with ${\bf b_1}$ (shorter) and ${\bf b_2}$ (longer,
or, for the first step only, at least as long as ${\bf b_1}$). Next, we reduce the longer vector ${\bf b_2}$ by adding or subtracting
multiples of the shorter vector ${\bf b_1}$, to get ${\bf b_3}={\bf b_2}-q{\bf b_1}$, where the integer $q$ is chosen to make $||{\bf b_3}||$
as short as possible. For ${\bf b_1}=(3,1)^{\sf T}$ and ${\bf b_2}=(5,3)^{\sf T}$, illustrated on the left of Figure~\ref{2step},
we get ${\bf b_3}=(-1,1)^{\sf T}$. Notice that $||{\bf b_3}||=\sqrt{2}$, so that ${\bf b_1}$ and ${\bf b_2}$ have exchanged places; the
shortened version ${\bf b_3}$ of ${\bf b_2}$ is now {\it strictly} shorter than ${\bf b_1}$. This means that we should exchange
the vectors, setting ${\bf b'_1=b_3}$ and ${\bf b'_2=b_1}$, and repeat. The second step of the algorithm is shown on the right of
Figure~\ref{2step}. The new shortened vector is ${\bf b'_3}=(2,2)^{\sf T}$, of length $\sqrt{8}$.

\noindent But now there is a problem. This time, ${\bf b'_1}$ and ${\bf b'_2}$ have {\it not} exchanged places; the adjusted longer vector
${\bf b'_3}$ is {\it still} longer than (or, in general, at least as long as) the shorter vector ${\bf b'_1}$. The resulting basis,
consisting of ${\bf b''_1=b'_1}$ and ${\bf b''_2=b'_3}$, is said to be {\it Lagrange-reduced}. This means that ${\bf b''_1}$ and
${\bf b''_2}$ satisfy
\[
||{\bf b''_1}||\le||{\bf b''_2}||\le||{\bf b''_2-qb''_1}||
\]
for all integers $q$. Lagrange's algorithm terminates here. (In this example, the final basis is orthogonal, but in general it need not be.)

\bigskip

\noindent It is nontrivial that Lagrange's algorithm does yield two basis vectors of minimal length, which are thus the closest and
second-closest lattice points to the origin. The proof, left as an exercise in~\cite{Gal}, and presented in full in~\cite{Bright}, is
about half a page long.

\bigskip

\noindent Lagrange's algorithm is very reminiscent of the Euclidean algorithm, with two quantities leapfrogging each other in a race
to the bottom. Consequently, it may now be clear to the reader that (as stated without proof in~\cite{Benne}) Brillhart's algorithm is
just Lagrange's algorithm applied to Grace's lattice~\eqref{eq:Grace}, starting with ${\bf b_1}=(u,1)^{\sf T}$ and ${\bf b_2}=(p,0)^{\sf T}$.
But the algorithms are in fact different, even for the case $p=73,u=27$ considered earlier. Before proceeding, the reader is encouraged to
run the two algorithms for this special case.

\section{The plot thickens}

\noindent The basic connection between Lagrange's algorithm and Brillhart's is that, {\it at least initially}, the $x$-components of
the sequence of vectors produced by Lagrange's algorithm (applied to $(p,0)^{\sf T}$ and $(u,1)^{\sf T}$) are the remainders in the Euclidean
algorithm (applied to $p$ and $u$). Indeed, if the $y$-components are small, we can neglect them, so that Lagrange's algorithm is
approximated by its one-dimensional version, which is precisely the Euclidean algorithm.

\bigskip

\begin{figure}[htp!]
\centering
\includegraphics[width=\linewidth]{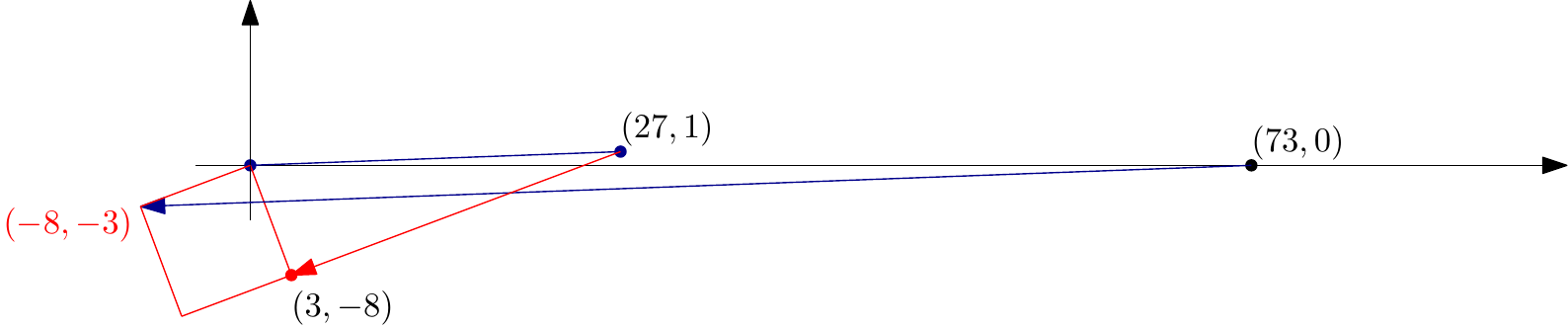}
\caption{\small Lagrange's algorithm}
\label{lagrange}
\end{figure}

\noindent Or is it? Figure~\ref{lagrange} shows the result of applying Lagrange's algorithm to the vectors $(73,0)^{\sf T}$ and $(27,1)^{\sf T}$.
It does indeed produce (essentially) the same answer as Brillhart's. But look at the $x$-components. Instead of the expected sequence
$73,27,19,8,3$, we see $73,27,3,-8$. The number 19 is absent entirely, and 8 and 3 are replaced by 3 and -8. So the two algorithms are different.
They are also different when, for example, $p=277$.

%

\bigskip

\noindent The resolution to this paradox is that the one dimensional version of Lagrange's algorithm is {\it not} the ordinary
version of Euclidean algorithm with non-negative remainders. Instead, it is the version where the remainders are chosen to minimize their
{\it absolute values}. So the obvious way to connect the algorithms is to modify Brillhart's algorithm to use this version of Euclid's
algorithm instead. Then it will run even faster.

\bigskip

\noindent The catch is that the fast Brillhart algorithm doesn't work, in the sense that it doesn't produce the required goods! Taking $p=277$
and $u=60$, we get
\begin{align*}
277&=5\cdot 60-23\\
60&=3\cdot 23-9\\
23&=3\cdot 9-4
\end{align*}
which misses the solution $9^2+14^2=277$ entirely.

\begin{figure}[htp!]
\centering
\includegraphics[width=\linewidth]{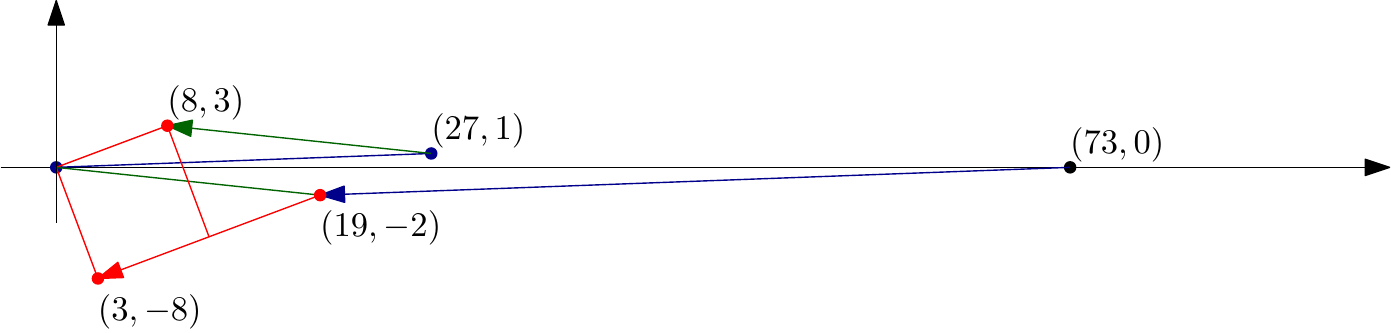}
\caption{\small Lagrange's algorithm, restricted to $x>0$}
\label{m-lagrange}
\end{figure}

\noindent Fortunately, there is another possibility. We can modify {\it Lagrange's} algorithm instead. And the obvious way to do this
is to insist that every vector it generates has positive $x$-component. Figure~\ref{m-lagrange} shows the result of applying this
modification of Lagrange's algorithm with initial vectors $(73,0)^{\sf T}$ and $(27,1)^{\sf T}$ This time, the sequence of $x$-components
does match the (positive) remainders. And, finally, it is not hard to modify the proof from~\cite{Bright} to show that the modified
Lagrange algorithm also terminates with a Lagrange-reduced basis, at least for the lattices under consideration.

\section{A calculation}

\noindent There is one final hurdle to clear. These approximations are only valid when the $y$-components of the relevant vectors are
small. Once they start growing, the shortest vector ${\bf b_3}$ from (the modified) Lagrange algorithm might not be the one with the
smallest $x$-component. We will show that, in fact, these vectors coincide, until both algorithms terminate with the same result.

\bigskip

\noindent To do this, it will be convenient to change notation slightly, and to work with points in $\R^2$, instead of their position
vectors. Let's introduce the notation by restating both algorithms.

\bigskip

\noindent {\bf Ingredients:} $p\equiv 1\pmod 4$ is a prime, and $u$ satisfies $u^2\equiv -1\pmod p$ with $1\le u<p/2$.

\bigskip

\noindent {\bf Initialization:} We start with ${\bf x}_0=(x_0,y_0)={\bf x}'_0=(p,0)$ and ${\bf x}_1=(x_1,y_1)={\bf x}'_1=(u,1)$.

\bigskip

\noindent {\bf Brillhart's Algorithm~\cite{Brill}:} For $i=1,2,\ldots$, we write $x_{i-1}=q_ix_i+r_i$ with $0\le r_i<x_i$, and set
\[
{\bf x}_{i+1}={\bf x}_{i-1}-q_i{\bf x}_{i}=(x_{i-1}-q_ix_i,y_{i-1}-q_iy_i)=(r_i,y_{i-1}-q_iy_i).
\]
Note that here we have extended Brillhart's algorithm to produce points in $\R^2$, and we have also extended it to run until Euclid's
algorithm terminates. The sequence $(x_i)$ is just the sequence of remainders $(r_i)$ from the Euclidean algorithm applied to $p$ and $u$,
while the sequence $(y_i)$ is also related to the Euclidean algorithm: we have $r_i=s_ip+y_iu$ for another sequence $(s_i)$ (see~\cite{Wagon}).

\bigskip

\noindent For example, with $p=193$ we have $u=81$, and successive ${\bf x}_i$ are given by
\[
(193,0),(81,1),(31,-2),(19,5),(12,-7),(7,12),(5,-19),(2,31),(1,-81),(0,193).
\]
The symmetry is related to the fact that $u^2\equiv -1\pmod p$; we will derive it as a consequence of the arguments below.

\bigskip

\noindent {\bf The Modified Lagrange (ML) Algorithm:} For $i=1,2,\ldots$, as long as $||{\bf x'}_i||\le ||{\bf x'}_{i-1}||$, write
\[
q'_i=\left\lfloor\frac{{\bf x'}_i\cdot {\bf x'}_{i-1}}{||{\bf x'}_i||^2}\right\rfloor,
\]
and set
\[
{\bf x'}_{i+1}={\bf x'}_{i-1}-q'_i{\bf x'}_{i}.
\]

\bigskip

\noindent Our aim is to show that, with the above choices of $p$ and $u$, the two sequences do indeed coincide, until the ML algorithm
terminates at ${\bf x}_t=(x_t,y_t)$. (In the example above, $(x_t,y_t)=(7,12)$.)
In other words, we need to show that $q_i=q'_i$, until the ML algorithm terminates.

\bigskip

\noindent We recall some details from earlier. First, the integer lattice $L$ given by
\[
L=\{(x,y)\in\Z^2: x\equiv uy\pmod p\}
\]
is a square lattice with determinant $p$, so that, for all $(x,y)\in L$, we have $p\mid(x^2+y^2)$,
and, if $(x,y)\in L$ is a point closest to the origin, then $p=x^2+y^2$. Consequently, if Lagrange's algorithm
terminates at ${\bf x}_t=(x_t,y_t)$, then $p=||{\bf x}_t||^2=x_t^2+y_t^2$.

\medskip

\noindent We now observe that each of the points ${\bf x}_i$ (and ${\bf x}'_i$) lies on $L$.
Consequently, for all $i\ge 1$,
\begin{equation}
p\mid(x_i^2+y_i^2).
\end{equation}
Moreover, it is easy to see by induction that, for all $i\ge 1$,
\begin{equation}
x_{i-1}y_i-x_iy_{i-1}=p,
\end{equation}
and so
\begin{align*}
(x_{i-1}x_i+y_{i-1}y_i)^2&=(x_{i-1}^2+y_{i-1}^2)(x_i^2+y_i^2)-(x_{i-1}y_i-x_iy_{i-1})^2\\
&=(x_{i-1}^2+y_{i-1}^2)(x_i^2+y_i^2)-p^2,
\end{align*}
so that
\begin{equation}
p\mid(x_{i-1}x_i+y_{i-1}y_i).
\end{equation}

\medskip

\noindent We can use these properties to establish the symmetry observed earlier. Suppose the (extended) Brillhart algorithm
terminates at $(x_T,y_T)$. Then, since the $x_i$ are just the reminders $r_i$ in the Euclidean algorithm, we must have $x_T=0$
and $x_{T-1}=1$. By (2) we get $y_i=p$, and, since $(x_{T-1},y_{T-1})\in L$, we then get $y_{i-1}=-u$. We can now establish
the required symmetry by induction (see also~\cite{Wagon}).

\medskip

\noindent Now, it is easy to show, by induction, that the $y_i$ alternate in sign, and that the sequence $|y_i|$ is increasing.
From the remarks above, the sequence $(|y_i|)$ is just the sequence $(x_i)$ in reverse order. Consequently, we can stop iterating
once $x_i\le |y_i|$. Thus, in the following, we will assume that $x_i>|y_i|$.

\medskip

\noindent{\bf Claim:} For all $i\ge 1$, $q_i=q'_i$.

\medskip

\noindent{\bf Proof:} We use induction on $i$. The induction starts. For the induction step, we need to show that
\[
q_i:=\left\lfloor\frac{x_{i-1}}{x_i}\right\rfloor=\left\lfloor\frac{x_{i-1}x_i+y_{i-1}y_i}{x^2_i+y_i^2}\right\rfloor=:q'_i.
\]
For simplicity of notation, write $a,b,c,d$ for $x_{i-1},y_{i-1},x_i,y_i$ respectively. Then we need that
\[
\left\lfloor\frac{a}{c}\right\rfloor=\left\lfloor\frac{ac+bd}{c^2+d^2}\right\rfloor.
\]
Suppose first that $q'_i<q_i$. Then there is an integer $m$ with
\begin{equation}
\frac{ac+bd}{c^2+d^2}<m\le\frac{a}{c}.
\end{equation}
Now, from the above, we know there are integers $K$ and $L$
\[
ac+bd=pK{\rm\ and\ }c^2+d^2=pL.
\]
Multiplying (4) by $c(c^2+d^2)$,
and noting that $y_i=d<c=x_i$, we get
\begin{align*}
cpK&=c(ac+bd)<mc(c^2+d^2)\\
&=mcpL\le a(c^2+d^2)\\
&=c(ac+bd)+d(ad-bc)\\
&=cpK+dp<cpK+cp\\
&=cp(K+1),
\end{align*}
so that, dividing both sides by $cp$, we have
\[
K<mL<K+1,
\]
a contradiction.

\medskip

\noindent The case where $q_i<q'_i$ is handled similarly; this corresponds to the case $y_i=d<0$.

\section{Exercise}

\noindent There is one remaining mystery. We've shown that Brillhart's algorithm is equivalent to the
modified Lagrange algorithm. The original Lagrange algorithm, meanwhile, corresponds to the modified
Brillhart algorithm, but the two can't quite be equivalent, since the former works and the latter doesn't.
Explain why!

\section{Acknowledgment}

\noindent Both authors are extremely grateful to the second author's teacher and mentor, Peter Shiu, whose
paper~\cite{Shiu} was the inspiration for this project, and whose subsequent communication and feedback
has been indispensable.


\begin{thebibliography}{99}

\bibitem{Bright} C. Bright,
Algorithms for lattice bass reduction , available at:

\noindent https://cs.uwaterloo.ca/~cbright/reports/latticealgs.pdf.

\bibitem{Brill} J. Brillhart,
Note on representing a prime as a sum of two squares,
{\it Math. Comp.} {\bf 26} (1972), 1011--1013.

\bibitem{Gal} S. Galbraith,
Mathematics and Public Key Cryptography,
Cambridge University Press, 2012.

\bibitem{Grace} J. Grace,
The four square theorem,
{\it J. Lond. Math. Soc.} {\bf 2} (1927), 3--8.

\bibitem{Shiu} P. Shiu,
The two squares theorem of Fermat for a prime $p\equiv 1\pmod 4$,
{\it Newsletter Lond. Math. Soc.} {\bf 494}, 26--31.

\bibitem{Wagon} S. Wagon,
The Euclidean algorithm strikes again,
{\it Amer. Math. Monthly} {\bf 97} (1990), 125--129.

\bibitem{Benne} B. de Weger,
Lagrange's algorithm strikes again, available at:

\noindent https://www.win.tue.nl/$\sim$bdeweger/downloads/sumof2squares.pdf.

\end{thebibliography}
\end{document}